\documentclass[10pt,a4paper]{article} 

\usepackage{a4wide}
\usepackage[latin1,utf8]{inputenc}
\usepackage[T1]{fontenc} 
\usepackage{textcomp} 
\usepackage[english]{babel}
\usepackage[autolanguage]{numprint} 
\usepackage{hyperref} 
\usepackage{graphicx} 
\usepackage{verbatim} 
\usepackage[all]{xy}
\usepackage{float}

\usepackage{lmodern}
\usepackage{mathrsfs}

\usepackage{amsmath,amssymb,amsthm} 
\usepackage{tikz,tkz-tab} 

\usepackage{mathdots}
\usepackage{makeidx}
\usepackage{stmaryrd} 

\usepackage{enumitem}
\setenumerate[0]{label=(\alph*)}

\renewcommand\epsilon\varespilon 
\renewcommand\phi\varphi 

\newcommand\bP{\mathbf{P}} 
\newcommand\bL{\mathbf{L}} 
\newcommand\ba{\mathbf{a}}  
\newcommand\bx{\mathbf{x}}  
\newcommand\by{\mathbf{y}}  
\newcommand\bz{\mathbf{z}}  
\newcommand\bu{\mathbf{u}}  
\newcommand\CCC{\mathcal{C}} 
\newcommand\CP{\textsl{P}} 
\newcommand\CR{\mathrm{R}} 
\newcommand\CX{\mathrm{X}}

\newcommand\CL{\textsl{L}}

\newcommand\ee{\varepsilon}
\newcommand\eee{\mu}

\newcommand\GrO{\mathcal{O}} 

\newcommand\hlambda{\widehat{\lambda}} 
\newcommand\hnu{\widehat{\nu}} 

\newcommand\ie{\textsl{i.e. }} 
\newcommand\kap{\kappa} 

\newcommand\NN{\mathbb{N}} 
\newcommand{\norm}[1]{\|#1\|} 
\newcommand\QQ{\mathbb{Q}} 
\newcommand\RR{\mathbb{R}} 

\newcommand\Sturm{\mathcal{S}} 

\newcommand\ZZ{\mathbb{Z}} 

\newcommand\pu{\underline{\psi}} 
\newcommand\po{\overline{\psi}} 

\newcommand\lambdaL{\hlambda_{\min}}
\newcommand\nuL{\hnu_{\min}}

\newcommand\NNN{\mathrm{N}}

\newcommand\VVV{\mathcal{V}}

\theoremstyle{definition} 
\newtheorem{Def}{Definition}[section]

\theoremstyle{plain} 
\newtheorem{Prop}[Def]{Proposition} 
\newtheorem{Lem}[Def]{Lemma} 
\newtheorem{Thm}{Theorem}[section] 

\theoremstyle{remark} 
\newtheorem{Rem}[Def]{Remark} 

\newenvironment*{Dem}{\noindent{\bf Proof}}{\hfill$\square$} 

\numberwithin{equation}{section} 

\newcommand{\Addresses}{{
  \bigskip
  \footnotesize
  \noindent A.~Poëls, \\ \textsc{
Département de Mathématiques, Université d'Ottawa, 150 Louis-Pasteur, Ottawa, Ontario K1N 6N5, Canada}
\par\nopagebreak
  \noindent \textit{E-mail}: \texttt{anthony.poels@uottawa.ca}
}}

\newcommand{\MSC}{{
  \footnotesize
  \textbf{MSC~2010}: 11J13(Primary), 11H06 (Secondary), 11J82.
}}

\newcommand{\keysW}{{
  \footnotesize
  \textbf{Keywords}: Diophantine approximation, parametric geometry of numbers, simultaneous approximation, exponents of Diophantine approximation
}}

\newcommand{\Ack}{{
  \footnotesize

  \textbf{Acknowledgements}: I am very grateful to Stéphane Fischler and Damien Roy for giving me a lot of
  feedback on this work.
}}

\title{A new exponent of simultaneous rational approximation}  
\author{Anthony Poëls} 
\date{} 

\begin{document} 
\maketitle

\begin{abstract}
\noindent We introduce a new exponent of simultaneous rational approximation $\lambdaL(\xi,\eta)$ for pairs of real numbers $\xi,\eta$, in complement to the classical exponents $\lambda(\xi,\eta)$ of best approximation, and $\hlambda(\xi,\eta)$ of uniform approximation. It generalizes Fischler's exponent $\beta_0(\xi)$ in the sense that $\lambdaL(\xi,\xi^2) = 1/\beta_0(\xi)$ whenever $\lambda(\xi,\xi^2) = 1$. Using parametric geometry of numbers, we provide a complete description of the set of values taken by $(\lambda,\lambdaL)$ at pairs $(\xi,\eta)$ with $1,\xi,\eta$ linearly independent over $\QQ$.

\end{abstract}

\MSC

\keysW

\Ack

\section{Introduction}

Let $\xi$ and $\eta$ be non-zero real numbers. The following simultaneous approximation problem has been intensively studied during the last decades:

\bigskip

\noindent\textbf{Problem $E_{\lambda,X}$:} Given $\lambda>0$ and $X\geq 1$, we search for solutions $(x_0,x_1,x_2)\in\ZZ^3\setminus\{0\}$ of the system
\begin{align*}
1\leq |x_0|\leq X\quad\textrm{and}\quad\max(|x_0\xi-x_1|,|x_0\eta-x_2|) \leq X^{-\lambda}.
\end{align*}
We denote by $\lambda(\xi,\eta)$ (resp. $\hlambda(\xi,\eta)$) the supremum of real numbers $\lambda$ for which $E_{\lambda,X}$ admits a non-zero integer solution for arbitrarily large values of $X$ (resp. for each sufficiently large value of $X$). For all real numbers $\xi,\eta$, we have
\[
    \lambda(\xi,\eta) \geq \hlambda(\xi,\eta) \geq \frac{1}{2},
\]
the right-hand side inequality following from Dirichlet's box principle (or, equivalently, Minkowski’s theorem). The study of such Diophantine exponents of approximation goes back to Jarn\'ik and Khinchin, see \cite{bugeaud2015exponents} for a well supplied account of the topic. In this paper, we consider the following variant:

\bigskip

\noindent\textbf{Problem $E_{\lambda,\eee,X}$:} Given $\lambda>0$, $\eee\geq 0$ and $X>1$, we search for solutions $(x_0,x_1,x_2)\in\ZZ^3\setminus\{0\}$ of the system
\begin{align*}
1\leq |x_0|\leq X\quad\textrm{and}\quad\max(|x_0\xi-x_1|,|x_0\eta-x_2|) \leq \min(X^{-\lambda},|x_0|^{-\eee}).
\end{align*}
This was introduced by Fischler in \cite{fischler2007palindromic} in the special case where $\eta=\xi^2$.
For $0\leq \eee < \lambda(\xi,\eta)$, we denote by $\hlambda_\eee(\xi,\eta)$ the supremum of the real numbers $\lambda$ for which $E_{\lambda,\eee,X}$ admits a non-zero integer solution for each sufficiently large value of $X$. 
Note that the map $\eee\mapsto\hlambda_\eee(\xi,\eta)$ is non-increasing. We define
\begin{align}
\label{Eq Def lambda_min}
    \lambdaL(\xi,\eta) = \inf_{0<\eee < \lambda(\xi,\eta)} \hlambda_\eee(\xi,\eta) = \lim_{\eee\rightarrow \lambda(\xi,\eta)^-}\hlambda_\eee(\xi,\eta).
\end{align}


\noindent See Remark \ref{Remarque: interpretation lambda_min} and \eqref{Eq formule nu epsilon comme lim inf avec log} for an interpretation of $\lambdaL$. Note that for $\eee=0$ we have $\hlambda_0(\xi,\eta) = \hlambda(\xi,\eta)$, so that
\[
    \lambdaL(\xi,\eta) \leq \hlambda(\xi,\eta).
\]
More generally, we have $\hlambda_\eee(\xi,\eta) = \hlambda(\xi,\eta)$ for any $\eee < \hlambda(\xi,\eta)$. In particular, if $\hlambda(\xi,\eta) = \lambda(\xi,\eta)$, Definition~\eqref{Eq Def lambda_min} gives $\lambdaL(\xi,\eta) = \hlambda(\xi,\eta) = \lambda(\xi,\eta)$. Yet, it is well-known that $\hlambda(\xi,\eta) = \lambda(\xi,\eta) = 1/2$ for almost all $(\xi,\eta)$ with respect to the Lebesgue measure on $\RR^2$ (see \cite[§2]{bugeaud2005exponents}). We thus have the following result:
\begin{Thm}
\label{Thm valeur générique lambda_min}
For almost all real numbers $\xi,\eta$ (with respect to the Lebesgue measure on $\RR^2$), we have
\[
    \lambdaL(\xi,\eta) = \frac{1}{2}.
\]
\end{Thm}


\noindent The goal of this paper is to give an interpretation of the exponents $\hlambda_\eee(\xi,\eta)$ and $\lambdaL(\xi,\eta)$ in the setting of parametric geometry of numbers and to prove the following description for the spectrum of the pair $(\lambda,\lambdaL)$, \ie the set of values taken by $(\lambda,\lambdaL)$ at pairs $(\xi,\eta)$ with $1,\xi,\eta$ linearly independent over $\QQ$.

\begin{Thm}
\label{Thm general spectre (lambda_min,lambda)}
For any $\xi,\eta\in\RR$ with $1,\xi,\eta$ linearly independent over $\QQ$, we have either $\lambdaL(\xi,\eta) = \lambda(\xi,\eta) = 1/2$, or
\begin{align}
\label{Eq Thm ppal 0}
 0 \leq \lambdaL(\xi,\eta) \leq 1,\quad \frac{1}{2} < \lambda(\xi,\eta) \leq +\infty \quad\textrm{and}\quad \frac{\lambdaL(\xi,\eta)^2}{1-\lambdaL(\xi,\eta)} \leq\lambda(\xi,\eta).
\end{align}
Conversely, for any $\hlambda\in\RR$ and any $\lambda\in\RR\cup\{+\infty\}$ satisfying either $\hlambda = \lambda = 1/2$, or
\begin{align}
\label{Eq Thm ppal}
    0 \leq \hlambda \leq 1,\quad \frac{1}{2} < \lambda \leq +\infty \quad\textrm{and}\quad \frac{\hlambda^2}{1-\hlambda} \leq\lambda,
\end{align}
there exist two real numbers $\xi$ and $\eta$, with $1,\xi,\eta$ linearly independent over $\QQ$, such that
\[
    \lambda(\xi,\eta) = \lambda\quad\textrm{and}\quad \lambdaL(\xi,\eta) = \hlambda.
\]
\end{Thm}

\noindent Laurent computed the spectrum of $(\lambda,\hlambda)$ in \cite{laurent2006exponents} (see Corollary~2 of \cite{laurent2006exponents}). He proved that for any $\xi,\eta$ with $1,\xi,\eta$ linearly independent over $\QQ$, we have
\begin{align}
\label{Eq relation Laurent}
     \frac{1}{2} \leq \hlambda(\xi,\eta) \leq 1, \quad \frac{\hlambda(\xi,\eta)^2}{1-\hlambda(\xi,\eta)} \leq \lambda(\xi,\eta) \leq +\infty,
\end{align}
and that \eqref{Eq relation Laurent} describe entirely the spectrum of $(\lambda,\hlambda)$. Since $\lambdaL\leq \hlambda$, the inequalities \eqref{Eq Thm ppal 0} are implied by \eqref{Eq relation Laurent} together with $\lambda(\xi,\eta) > 1/2$. It would be interesting to study the joint spectrum of $(\lambda,\hlambda,\lambdaL)$.\\

In \cite{fischler2007palindromic} Fischler introduced a new exponent of approximation $\beta_0(\xi)$ for each real number $\xi$. When $\lambda(\xi,\xi^2)<1$, he defined $\beta_0(\xi) = +\infty$. Otherwise he set $\beta_0(\xi) = \lim_{\ee\rightarrow 0^+}\beta_\ee(\xi)$, with $\beta_\ee(\xi) = 1/\hlambda_{1-\ee}(\xi,\xi^2)$ (for $0<\ee \leq 1$). Then he studied in depth the real numbers $\xi$ for which $\beta_0(\xi) <2$. For those numbers, the exponent $\beta_0(\xi)$ and $\lambdaL(\xi,\xi^2)$ are related as follows.

\begin{Lem}
If $\beta_0(\xi) < 2$, then $\lambda(\xi,\xi^2) = 1$ and $\beta_0(\xi) = 1/\lambdaL(\xi,\xi^2)$.
\end{Lem}

\begin{Dem}
Let $\xi$ be such that $\beta_0(\xi) < 2$. Then we have $\lambda(\xi,\xi^2) \geq 1$. In general, the inequality $1/\beta_0(\xi) \leq \hlambda(\xi,\xi^2)$ holds, so that $\hlambda(\xi,\xi^2) > 1/2$. This implies that $\lambda(\xi,\xi^2) \leq 1$. This result can be obtained from Davenport and Schmidt's work by generalizing Lemmas~2 and~6 of \cite{davenport1969approximation} (see for example \cite[Corollaire 6.2.7]{poelsPhD}); it is also a corollary of a recent result due to Schleischitz \cite[Theorem 1.6]{schleischitz2016spectrum}. Finally, this shows that $\lambda(\xi,\xi^2) = 1$. In this case, we have $\beta_0(\xi) = 1/\lambdaL(\xi,\xi^2)$ by definition of $\beta_\ee(\xi)$ ($0<\ee<1$).

\end{Dem}

Let $\VVV$ denotes the set
\begin{equation}
\label{def ensemble VVV}
    \VVV = \big\{(\xi,\eta)\;|\; \textrm{$1,\xi,\eta$ linearly independent over $\QQ$ and $\lambda(\xi,\eta) = 1$} \big\}.
\end{equation}
Applying  Theorem~\ref{Thm general spectre (lambda_min,lambda)} with  $\lambda=1$, we obtain the following result.

\begin{Thm}
\label{Thm spec beta(xi,eta)}
    With the above notation, the set of values taken by $1/\lambdaL$ at pairs $(\xi,\eta)\in \VVV$ is $[\gamma,+\infty]$, where $\gamma = (1+\sqrt 5)/2$ denotes the golden ratio.
\end{Thm}

\noindent The situation is radically different for the pairs $(\xi,\xi^2)$. Following \cite{bugeaud2005exponentsSturmian} let us denote by $\Sturm$ the set of all values $\displaystyle \sigma= 1/\limsup_{k\rightarrow+\infty}[s_{k+1};s_k,\dots,s_1]$ where $(s_k)_{k\geq1}$ runs through all sequences of positive integers (here $[a_0 ; a_1,a_2,\dots]$ denotes the continued fraction whose partial quotients are $a_0, a_1,\dots$). The largest element of $\Sturm$ is $\frac{1}{\gamma}$. The values immediately below have been described by Cassaigne \cite{cassaigne1999limit}. They constitute a decreasing sequence of quadratic numbers converging to the largest accumulation point $s\approx 0.3867\dots$ of $\Sturm$. Also note that Cassaigne has shown in \cite{cassaigne1999limit} that this set has empty interior. Elements of $\Sturm$ appear in the description of the classical exponents of approximation to Sturmian continued fractions, studied by Bugeaud and Laurent in \cite{bugeaud2005exponentsSturmian}, and of Sturmian type numbers (see \cite{poels2017exponents}). The set $\Sturm$ is related to the spectrum of $\beta_0$ by the following result (see \cite{fischler2007palindromic}):

\begin{Thm}[Fischler, 2007]
\label{Thm spectre beta_0 < sqrt 3}
Let us set $\Sturm_0 = \{\beta_0(\xi)\;|\; \textrm{$\xi\in\RR$ not algebraic of degree $\leq 2$}\}$. Then we have
\[
    \Sturm_0\cap[\gamma,\sqrt 3) = \Big\{1+\frac{1}{1+\sigma}\;|\; \sigma\in\Sturm\Big\}\cap(1,\sqrt 3).
\]
\end{Thm}

\noindent In view of the description of $\Sturm$ given above, the smallest element of $\Sturm_0$ in $[\gamma,+\infty)$ is therefore $\gamma$ and the values immediately above constitute an increasing sequence of quadratic numbers converging to the smallest accumulation point \mbox{$1.721\dots < \sqrt 3$} of $\Sturm_0$. Thus, Theorem~\ref{Thm spec beta(xi,eta)} implies that $\{1/\lambdaL(\xi,\eta)\;|\;(\xi,\eta)\in\VVV\}\cap[\gamma,\sqrt 3]$ is the full interval $[\gamma,\sqrt 3]$ (where $\VVV$ is defined by \eqref{def ensemble VVV}), whereas
\[
    \{\beta_0(\xi)\;|\;\textrm{$\xi\in\RR$ not algebraic of degree $\leq 2$}\}\cap[\gamma,\sqrt 3] = \{1/\lambdaL(\xi,\xi^2)\;|\;(\xi,\xi^2)\in\VVV\}\cap[\gamma,\sqrt 3]
\]
has empty interior and its complement in $[\gamma,\sqrt 3]$ has non-empty interior by Theorem~\ref{Thm spectre beta_0 < sqrt 3}.\\


\noindent In Section~\ref{section exposants lambda epsilon} we study the ``rigidity'' of the exponents $\hlambda_\eee$ and we recall the notion of \textsl{minimal points}, which is useful to compute the exponents. We use parametric geometry of numbers to prove Theorem~\ref{Thm general spectre (lambda_min,lambda)}. In section~\ref{section PGN} we briefly recall the elements of that theory and we provide a parametric version of the exponent $\lambdaL(\xi,\eta)$. The proof of Theorem~\ref{Thm general spectre (lambda_min,lambda)} is given in Section~\ref{section proof thm principal}. Section~\ref{section: finals remarks and open questions} is devoted to some open questions about $\lambdaL$.

\section{Exponents $\hlambda_\eee$}
\label{section exposants lambda epsilon}

Unlike the classical exponents $\lambda(\xi,\eta)$ and $\hlambda(\xi,\eta)$, the exponent $\hlambda_\eee(\xi,\eta)$ may change if we perturbate the problem $E_{\lambda,\eee,X}$ slightly, for example by using $\norm{\bx}$ instead of $|x_0|$ (where $\norm{\cdot}$ is a fixed norm on $\RR^3$). However, as kindly pointed out to the author by Damien Roy, this happens only at the points $\eee$ at which the non-increasing map $\eee\mapsto \hlambda_\eee(\xi,\eta)$ is not continuous; this set of points is therefore countable. We formalize this claim in Proposition~\ref{Lemme lambda epsion bien def} below. Thus the exponent $\lambdaL(\xi,\eta)$ can be defined using any norm $\norm{\bx}$ instead of $|x_0|$ in $E_{\lambda,\eee,X}$.\\

\noindent Let $\xi,\eta\in\RR$ be two real numbers with $1,\xi,\eta$ linearly independent over $\QQ$. The exponents  $\hlambda_\eee(\xi,\eta)$ are defined as in the introduction. We denote by $\norm{\cdot}$ the usual Euclidean norm in $\RR^3$. If $f,g:I\rightarrow[0,+\infty)$ are two fonctions on a set $I$, we write $f\ll g$ (resp. $f \gg g$) to mean that there is a positive constant $c$ such that $f(x)\leq cg(x)$ (resp. $f(x) \geq cg(x)$) for each $x\in I$. We write $f\asymp g$ if both $f\ll g$ and $g \ll f$ hold. Let $\Delta,\NNN : \RR^3\rightarrow[0,+\infty)$ such that for any $\bx=(x_0,x_1,x_2)$
\[
    \Delta(\bx) \asymp \max\big(|x_0\xi-x_1|,|x_0\eta-x_2|\big)
\]
and
\[
    \NNN(\bx) \asymp \norm{\bx}\quad\textrm{if $\max\big(|x_0\xi-x_1|,|x_0\eta-x_2|\big) < 1$},
\]
(the implicit constants depend only on $\Delta,\NNN, \xi$ and $\eta$). Note that we may take $\NNN(x) = |x_0|$, although $\NNN$ is not a norm in this case. For $0\leq \eee< \lambda(\xi,\eta)$, we denote by $\hnu_\eee(\xi,\eta)$ the supremum of the real numbers $\nu$ for which the system
\begin{align}
\label{Eq system nu epsilon}
    \NNN(\bx)\leq X\quad \textrm{and}\quad \Delta(\bx) \leq \min\big(X^{-\nu},\NNN(\bx)^{-\eee}\big)
\end{align}
admits a non-zero integer solution for each sufficiently large value of $X$.
The map $\eee\mapsto\hnu_\eee(\xi,\eta)$ is non-increasing. We set
\[
    \nuL(\xi,\eta) = \inf_{0<\eee < \lambda(\xi,\eta)} \hnu_\eee(\xi,\eta) = \lim_{\eee\rightarrow \lambda(\xi,\eta)^-}\hnu_\eee(\xi,\eta).
\]

We have the following result:

\begin{Prop}
\label{Lemme lambda epsion bien def}
The non-increasing maps $\eee\mapsto\hnu_\eee(\xi,\eta)$ and $\eee\mapsto\hlambda_\eee(\xi,\eta)$  have the same set of discontinuities on $[0,+\infty)$ and they coincide outside of this set. Moreover we have:
\[
    \nuL(\xi,\eta) = \lambdaL(\xi,\eta).
\]
\end{Prop}

\begin{Dem}
Let us prove that $\hlambda_{\eee'}(\xi,\eta) \geq \hnu_{\eee}(\xi,\eta)$ for any $0\leq \eee' < \eee$. If $\hnu_\eee(\xi,\eta) = 0$ it is trivial. Now, suppose $\hnu_\eee(\xi,\eta) > 0$ and let $0<\lambda' < \lambda < \hnu_{\eee}(\xi,\eta)$. If $X$ is large enough, then  \eqref{Eq system nu epsilon} has a non-zero integer solution $\bx$ and this point $\bx$ is also solution of the problem $E_{\lambda',\eee',X}$ stated in the introduction. By letting $\lambda'$ tend to $\lambda$, then by letting $\lambda$ tend to $\hnu_{\eee}(\xi,\eta)$, it follows that $\hlambda_{\eee'}(\xi,\eta) \geq \hnu_{\eee}(\xi,\eta)$. Conversely, we also have $\hnu_{\eee'}(\xi,\eta) \geq \hlambda_{\eee}(\xi,\eta)$. In summary, we have shown that for any $\eee_1<\eee<\eee_2$, we have $\hlambda_{\eee_2}(\xi,\eta) \leq \hnu_\eee(\xi,\eta) \leq \hlambda_{\eee_1}(\xi,\eta)$, which yields $\hnu_\eee(\xi,\eta) = \hlambda_{\eee}(\xi,\eta)$ at each point where $\eee\mapsto\hlambda_\eee(\xi,\eta)$ (or $\eee\mapsto\hnu_\eee(\xi,\eta)$) is continuous.
\end{Dem}

\bigskip

To compute the exponent $\hnu_\eee(\xi,\eta)$ it is sufficient to consider only ``the best'' solutions of \eqref{Eq system nu epsilon}. Following Davenport and Schmidt \cite{davenport1967approximation}, \cite{davenport1969approximation}, we call a sequence of \textsl{minimal points} (associated to $\NNN$ and $\Delta$) a sequence of non-zero integer points $(\bx_i)_{i\geq 0}$ which satisfies
\begin{itemize}
\item[$\bullet$] $\NNN(\bx_1)<\NNN(\bx_2)<\dots$ and $\Delta(\bx_1) > \Delta(\bx_2)>\dots$,
\item[$\bullet$] For each non-zero $\bz\in\ZZ^3$, if $\NNN(\bz) < \NNN(\bx_{i+1})$, then $\Delta(\bz) \geq \Delta(\bx_i)$.
\end{itemize}
For simplicity, let us write $\CX_i = \NNN(\bx_i)$ and $\Delta_i = \Delta(\bx_i)$. Let $\eee\geq 0$, $\lambda>0$ and $X>0$, and suppose that $\bx\in\ZZ^3$ satisfies \eqref{Eq system nu epsilon}. If $\NNN(\bx) \gg 1$, then there is an index $i$ such that $\CX_i \leq \NNN(\bx) < \CX_{i+1}$. Since $\Delta_i \leq \Delta(\bx)$ and $\eee\geq 0$, the point $\bx_i$ is also solution of \eqref{Eq system nu epsilon}. Hence, for each $\eee < \lambda(\xi,\eta)$, the exponent $\hnu_\eee(\xi,\eta)$ is the supremum of the real numbers $\lambda$ such that for each $X$ large enough, there exists $i\geq 1$ for which
\[
    \CX_i \leq X\quad \textrm{and}\quad \Delta_i \leq \min\big(X^{-\lambda},\CX_i^{-\eee}\big).
\]
Let $0<i_1<i_2<\dots$ denote the sequence of indices $i$ such that $\Delta_i \leq \CX_i^{-\eee}$. Then
\begin{align}
\label{Eq formule nu epsilon comme lim inf avec log}
    \hnu_\eee(\xi,\eta) = \liminf_{k\rightarrow\infty} \frac{-\log(\Delta_{i_k})}{\log(\CX_{i_{k+1}})}.
\end{align}
For $\mu = 0$, we simply have
\begin{align*}
    \hlambda(\xi,\eta) = \hnu_0(\xi,\eta) = \liminf_{i\rightarrow\infty} \frac{-\log(\Delta_{i})}{\log(\CX_{i+1})}.
\end{align*}

\begin{Rem}
\label{Remarque: interpretation lambda_min}
Formula \eqref{Eq formule nu epsilon comme lim inf avec log} is similar to $(11)$ of \cite{fischler2007palindromic}. Roughly speaking, $\lambdaL(\xi,\eta)$ corresponds to $\hlambda(\xi,\eta)$ when we only take in account the exceptionally precise approximants, \ie solutions $\bx=(x_0,x_1,x_2)\in\ZZ^3$ of $E_{\lambda,X}$ with $\max(|x_0\xi-x_1|,|x_0\eta-x_2|)$ very close to $|x_0|^{-\lambda(\xi,\eta)}$; the quantity $1/\lambdaL(\xi,\eta)$ measures the ``maximal gap'' between two successive such very good approximants. If $\lambdaL(\xi,\eta) = 0$, then the ``maximal gap'' is big: for any $\ee>0$ small enough, there are infinitely many consecutive exceptionally precise approximants $\by$ and $\bz$, such that $\norm{\by} \leq \norm{\bz}^\ee$. If $\lambdaL(\xi,\eta) = 1$, then the ``maximal gap'' is small: for any consecutive exceptionally precise approximants $\by$ and $\bz$, $\norm{\by}$ is very close to $\norm{\bz}^{1/\lambda(\xi,\eta)}$.

\end{Rem}

\section{Parametric geometry of numbers}
\label{section PGN}

\subsection{The setting}
\label{subsection setting}

In this section we quickly recall the basics of the parametric geometry of numbers following Schmidt and Summerer \cite{Schmidt2009}, \cite{Schmidt2013} and Roy \cite{Roy_juin}. We use the setting of \cite{Roy_juin}. We denote by $\bx\wedge\by$ the standard vector product of two vectors $\bx,\by\in\RR^3$, by $\bx\cdot\by$ their standard inner product and by $\norm{\bx}$ the Euclidean norm of $\bx$. Fix $\bu\in\RR^3\setminus\{0\}$. For each $q\geq 0$ we set
\[
    \CCC_\bu(q)~:= \{\bx\in\RR^3\;;\; \norm{\bx} \leq 1, |\bx\cdot\bu|\leq e^{-q} \}\quad \textrm{and}\quad\mathcal{C}_\bu^*(q) := \{\bx\in\RR^3\;;\; \norm{\bx} \leq e^{q}, \norm{\bx\wedge\bu}\leq 1 \}.
\]
For $j=1,2,3$ we define a function $\CL_j:[0,+\infty)\rightarrow\RR$ by $\CL_{j}(q) = \log (\lambda_{j,\bu}(q))$, where $\lambda_{j,\bu}(q)$ denotes the $j$-th successive minimum of the convex body $\CCC_{\bu}(q)$ with respect to the lattice $\ZZ^3$. We set $\bL_\bu=(\CL_1,\CL_2,\CL_3)$. The functions $\CL_{j}$ are continuous, piecewise linear with slopes $0$ and $1$, and by Minkowski's second theorem they satisfy $\CL_1(q)+\CL_2(q)+\CL_3(q) = q + \GrO(1)$ (for any $q\geq 0$). For each $\bx\in\RR^3$, we further define $\lambda_\bx(q,\CCC_\bu(q))$ to be the smallest real number $\lambda\geq 0$ such that $\bx\in\lambda\CCC_\bu(q)$. When $\bx\neq 0$, this number is positive and so we obtain a function $\CL_\bx : [0,+\infty) \rightarrow\RR$ by putting $\CL_\bx(q) := \log(\lambda_\bx(q,\CCC_\bu(q)))$.
\noindent For $j=1,2,3$ we set
\begin{align*}
\po_j(\bu) = \po_j = \limsup_{q\rightarrow\infty}\frac{\CL_{j}(q)}{q}\quad \textrm{and}\quad \pu_j(\bu) = \pu_{j} = \liminf_{q\rightarrow\infty}\frac{\CL_{j}(q)}{q}.
\end{align*}
Similarly we define the function $\bL_\bu^*=(\CL_1^*,\CL_2^*,\CL_3^*)$, $\po_j^*$, $\CL_\bx^*$ ($\bx\in\RR^3\setminus\{0\}$), $\pu_j^*$ associated to the family of convex bodies $\CCC_\bu^*(q)$. For any non-zero $\bx\in\RR^3$ we have
\begin{align}
\label{Eq formule pour L_bx^*}
    \CL_\bx^*(q) = \max\big(\log \norm{\bx\wedge\bu},\log \norm{\bx}-q \big)\quad (q \geq 0).
\end{align}
\noindent The dual functions $\CL_j^*$ are related to the functions $\CL_j$ by Mahler's duality:

\begin{Prop} [Mahler]
\label{Prop Mahler}
For $j=1,2,3$ we have $\mathrm{L}_{j}(q) = -\mathrm{L}_{4-j}^*(q)+\GrO(1)$ for all $q>0$.
\end{Prop}
\noindent Thus
\begin{align}
\label{Eq dualité exposants}
    \pu_j= - \po_{4-j}^*\quad\textrm{and}\quad \po_j= - \pu_{4-j}^*\quad (j=1,2,3).
\end{align}

\noindent The following definition is that of a $3$-system (see \cite[Definition 4.1]{Roy_octobre}; this is an analog of a $(3,0)$-system for Schmidt and Summerer \cite{Schmidt2013}).

\begin{Def}
\label{Def n-système}
Fix a real number $q_0\geq0$. A $3$-system on  $[q_0,+\infty)$ is a continuous piecewise linear map $\bP=(\CP_1,\CP_2,\CP_3):[q_0,+\infty)\rightarrow\RR^3$  with the following properties:
\begin{enumerate}
\item \label{Def n-système condition 1} For each $q\geq q_0$, we have $0\leq \CP_1(q)\leq\CP_2(q)\leq \CP_{3}(q)$ and $\CP_1(q)+\CP_2(q)+\CP_3(q) = q$.
\item \label{Def n-système condition 2} If $H$ is a non-empty open subinterval of $[q_0,+\infty)$ on which $\bP$ is differentiable, then there is an integer $r$ ($1\leq r\leq 3$), such that $\CP_r$ has slope $1$ on $H$ while the other components $\CP_j$ of $\bP$ ($j\neq r$) are constant on $H$.
\item \label{Def n-système condition 3} If $q>q_0$ is a point at which $\bP$ is not differentiable and if the integers $r$ and $s$, for which $\CP_r$ has slope $1$ on $(q-\ee,q)$ and $\CP_s$ has slope $1$ on $(q,q+\ee)$ (for  $\ee>0$ small enough), satisfy $r<s$, then we have $\CP_r(q)=\CP_{r+1}(q)=\dots=\CP_s(q)$.
\end{enumerate}
\end{Def}

The following fondamental result was proved by Roy in \cite{Roy_juin}.

\begin{Thm}[Roy, $2015$]
\label{Thm Roy conjecture S&S}
For each non-zero point $\bu\in\RR^3$, there exist $q_0>0$ and a $3$-system $\bP$ on $[q_0,+\infty)$ such that $\norm{\bL_{\bu}-\bP}_{\infty}$ is bounded on $[q_0,+\infty)$. Conversely, for each  $3$-system $\bP$ on an interval $[q_0,+\infty)$, there exists a non-zero point $\bu\in\RR^3$ such that $\norm{\bL_{\bu}-\bP}_{\infty}$ is bounded on $[q_0,+\infty)$.
\end{Thm}

Following \cite[§3]{Schmidt2013} we define the \emph{combined graph} of a set of real valued functions defined on an interval $I$ to be the union of their graphs in $I\times\RR$. For a map $\bP:[c,+\infty)\rightarrow \RR^3$ and an interval $I\subset [c,+\infty)$, we also define the \emph{combined graph of $\bP$ on $I$} to be the combined graph of its components $P_1,P_2,P_3$ restricted to $I$.\\


We recall the following relationship between classical and parametric exponents (see \cite{Roy_juin}). For any $\bu=(1,\xi,\eta)$ with $\QQ$-linearly independent coordinates, we have
\begin{align}
\label{eq dico}
\big(\pu_3(\bu),\po_3(\bu)\big) = \Big(\frac{\hlambda(\xi,\eta)}{1+\hlambda(\xi,\eta)}, \frac{\lambda(\xi,\eta)}{1+\lambda(\xi,\eta)} \Big).
\end{align}

\subsection{Parametric formulation of $\lambdaL$}

\begin{Def}
\label{Def kappa(P)}
Let $c\geq 0$ and let $\CP:[c,+\infty)\rightarrow[0,+\infty)$ be an unbounded continuous piecewise linear function, with slopes $0$ and $1$. Let $(q_i)_{i\geq 0}$ be the increasing sequence of abscissas at which $\CP$ changes slope from $1$ to $0$. We suppose $(q_i)_{i\geq 0}$ infinite and define $\po(\CP)$, $\pu(\CP)$ by
\[
    \po(\CP) = \limsup_{q\rightarrow+\infty} \frac{\CP(q)}{q} = \limsup_{k\rightarrow\infty}\frac{\CP(q_k)}{q_k} \quad \textrm{and}\quad \pu(\CP) = \liminf_{q\rightarrow+\infty} \frac{\CP(q)}{q}.
\]
For each $\alpha < \po(\CP)$, let $(q_{i,\alpha})_{i\geq 0}$ be the (increasing) subsequence of all abscissas $q_k$ satisfying $q_k^{-1}\CP(q_k) \geq \alpha$. For each $i\geq 0$ we denote by $r_{i,\alpha}$ the abscissa of the intersection point of the horizontal line passing through $(q_{i,\alpha},\CP(q_{i,\alpha}))$ and of the line with slope $1$ passing through $(q_{i+1,\alpha},\CP(q_{i+1,\alpha}))$. We set
\[
    \kap_\alpha(\CP) = \liminf_{i\rightarrow+\infty}\frac{\CP(q_{i,\alpha})}{r_{i,\alpha}} \quad \textrm{and}\quad \kap(\CP) = \lim_{\alpha\rightarrow \po(\CP)} \kap_\alpha(\CP).
\]
Let $\CP^*:[c,+\infty)\rightarrow(-\infty,0]$ be an unbounded continuous piecewise linear function, with slopes $0$ and $-1$ and which changes from slope $-1$ to $0$ infinitely many times. In a dual manner, for $\displaystyle \alpha > \liminf_{q\rightarrow\infty} \CP^*(q)/q$ we define
\[
    \kap^*_\alpha(\CP^*) = -\kap_{-\alpha}(-\CP^*)\quad\textrm{and}\quad \kap^*(\CP^*) = -\kap(-\CP^*).
\]
\end{Def}

\noindent Note that $\kap(\CP)\leq \pu(\CP)$.

\begin{Lem}
\label{Lemma kappa(P) = kappa(R)}
Let $c\geq 0$ and let $\CP$, $\CR$ be two unbounded non-negative continuous piecewise linear functions defined on $[c,+\infty)$, with slopes $0$ and $1$, and which change from slope $1$ to $0$ infinitely many times. Suppose that $|\CP(q)-\CR(q)| = o(q)$ as $q$ tends to infinity. Then, the non-increasing maps $\alpha\mapsto \kap_\alpha(\CP)$ and $\alpha\mapsto \kap_\alpha(\CR)$ have the same set of discontinuities on $[0,1[$ and they coincide outside of this set. Moreover, we have $\kap(\CP) = \kap(\CR)$.
\end{Lem}

\begin{Dem}
Let $\CP$ and $\CR$ be as above. By hypothesis we have $\po(\CP) = \po(\CR)=:\po$ and $\pu(\CP) = \pu(\CR)=:\pu$. Fix $\alpha<\beta<\po$ and let us denote by $(q_i^\CP)_i$, $(q_{i,\beta}^\CP)_i$, $(r_{i,\beta}^\CP)_i$ (resp. $(q_i^R)_i$, $(q_{i,\alpha}^R)_i$, $(r_{i,\alpha}^R)_i$) the quantities associated by Definition~\ref{Def kappa(P)} to $\kap_\beta(\CP)$ (resp. $\kap_\alpha(\CR)$). Let us first prove that
\begin{align}
\label{Eq kappa essentiellement décroissante}
    \kap_\beta(\CP) \leq \kap_\alpha(\CR).
\end{align}

\noindent  Let $\ee>0$ be such that $\alpha+\ee < \beta$ 
and fix $i$ arbitrarily large. If $\CR(q)\geq \alpha q$ for each $q\in K_i:=[q_{i,\alpha}^\CR,q_{i+1,\alpha}^\CR]$, then we set $r = s = r_{i,\alpha}^\CR$. Otherwise $[r,s]$ denotes the maximal subinterval of $K_i$ on which $R(q)\leq \alpha q$. Let us write $A_1=(r,\CR(r))$, $A_2=(s,\CR(s))$ and $A_3=(r_{i,\alpha}^\CR, \CR(q_{i,\alpha}^\CR))$. The graph of $\CR$ above the interval $K_i$ is contained inside the triangle $(A_1A_2A_3)$. Let us denote by $\mathcal{D}_1$ (resp. $\mathcal{D}_2$) the horizontal line passing through the point $(r,\CR(r)+\ee r)$ (resp. the line with slope $1$ passing through $(s,\CR(s)+\ee s)$ (see Figure~\ref{figure graphe R et P calcul kappa})).
\begin{figure}[H]
    \begin{tabular}{l|r}
    \begin{tikzpicture}[scale=0.5]
        \clip(0.5,-1.5) rectangle (13.8,9);

        \draw [white, top color=gray!50,bottom color=white, shading angle = 45] (4.5,0)-- (8,3.5)-- (13,3.5)-- (13,0) -- cycle; 


        \draw [loosely dashed] (0,3.5) -- (13,3.5); 
        \draw (13, 3.5) node [above] {$\mathcal{D}_1$};
        \draw [loosely dashed] (4.5,0) -- ++(8.5,8.5); 
        \draw (12.8, 8.5) node [left] {$\mathcal{D}_2$};
        \draw (8,3.5) node {$\bullet$} node [above left] {$M$};

        \draw (12.5,5+1/3) node {slope $\alpha$}; 
        \draw (0,1.2) -- ++(13,13*0.4); 

        \draw [line width=2pt](0,2)-- (1,3)-- (4.5,3);
        \draw [dotted] (1,3) -- (1,0) node [below] {$q_{i,\alpha}^\CR$};

        \draw [densely dotted] (4.5,3) node {$\bullet$} node [below] {$A_1$} -- (8,3) node {$\bullet$} node [below] {$A_3$}-- (10+1/3,5+1/3) node {$\bullet$} node [below] {$A_2$};
        \draw [dotted] (4.5,3) -- (4.5,0) node [below] {$r$};
        \draw [dotted] (8,3) -- (8,0) node [below] {$r_{i,\alpha}^\CR$};
        \draw [dotted] (10+1/3,5+1/3) -- (10+1/3,0) node [below] {$s$};

        \draw [line width=2pt] (10+1/3,5+1/3)-- (12,7)-- (13,7) node [above] {$\CR$};
        \draw [dotted] (12,7) -- (12,0) node [below] {$q_{i+1,\alpha}^\CR$};
    \end{tikzpicture}

    &

    \begin{tikzpicture}[scale=0.5]
        \clip(0.5,-1.5) rectangle (13.8,9);

        \draw [white, top color=gray!50,bottom color=white, shading angle = 45] (4.5,0)-- (8,3.5)-- (13,3.5)-- (13,0) -- cycle; 


        \draw [loosely dashed] (0,3.5) -- (13,3.5); 
        \draw (13, 3.5) node [above] {$\mathcal{D}_1$};
        \draw [loosely dashed] (4.5,0) -- ++(8.5,8.5); 
        \draw (12.8, 8.5) node [left] {$\mathcal{D}_2$};
        \draw (8,3.5) node {$\bullet$} node [above left] {$M$};

        \draw (12.5,5+1/3-2.6) node {slope $\alpha$}; 
        \draw (2,2*0.4+1.2-2.6) -- ++(9.5,9.5*0.4); 

        \draw [line width=2pt](0,2)-- (1,3)-- (4.5,3);
        \draw [dotted] (1,3) -- (1,0) node [below] {$q_{i,\alpha}^\CR$};

        \draw [line width=2pt] (4.5,3) -- (8,3) node {$\bullet$} node [below] {$A_3$}-- (10+1/3,5+1/3);
        \draw [dotted] (8,3) -- (8,0) node [below] {$r = s = r_{i,\alpha}^\CR$};

        \draw [line width=2pt] (10+1/3,5+1/3)-- (12,7)-- (13,7) node [above] {$\CR$};
        \draw [dotted] (12,7) -- (12,0) node [below] {$q_{i+1,\alpha}^\CR$};
    \end{tikzpicture}
    \end{tabular}
    \caption{ \label{figure graphe R et P calcul kappa}
    Graph of $\CR$ on $[q_{i,\alpha}^\CR,q_{i+1,\alpha}^\CR]$}
\end{figure}
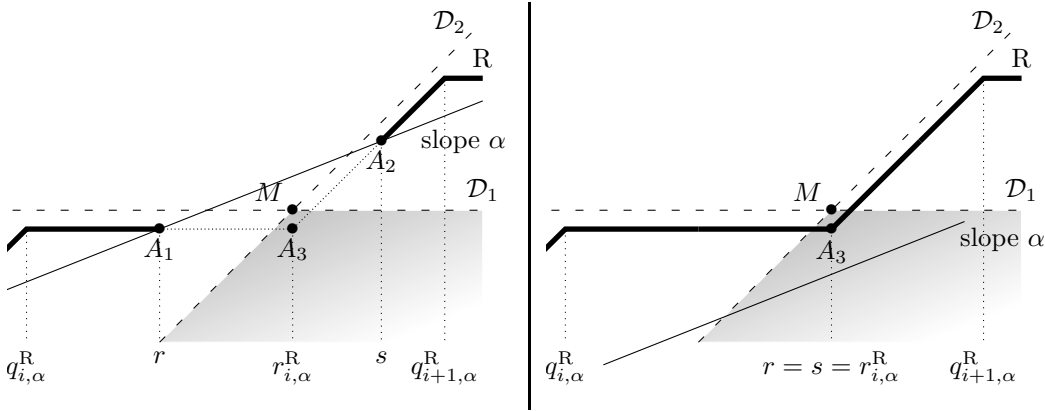

\noindent Now, let us define $j$ as the maximal index such that $q_{j,\beta}^\CP\leq r$. The horizontal line passing through $(q_{j,\beta}^P,\CP(q_{j,\beta}^P))$ lies below the line $\mathcal{D}_1$. If $r=s$, then $q_{j+1,\beta}^\CP \geq s$. Otherwise, if $i$ is large enough, then for each $q\in[r,s]$ we have:
\[
    \CP(q) \leq \CR(q) + \ee q \leq (\alpha+\ee) q  < \beta q,
\]
which implies that $q_{j+1,\beta}^\CP\geq s$. It follows that the line with slope $1$ passing through $(q_{j+1,\beta}^P,\CP(q_{j+1,\beta}^P))$ is below the line $\mathcal{D}_2$. As a consequence, the point $(r_{j,\beta}^\CP,\CP(q_{j,\beta}^\CP))$ lies in the area below $\mathcal{D}_1$ and $\mathcal{D}_2$ (the gray area of Figure~\ref{figure graphe R et P calcul kappa}). Let $M=(x_M,y_M)$ denote the intersection point of $\mathcal{D}_1$ and $\mathcal{D}_2$. Then we have
\[
    \frac{\CP(q_{j,\beta}^\CP)}{r_{j,\beta}^\CP} \leq \frac{y_M}{x_M} = \frac{\CR(q_{i,\alpha}^\CR)+\ee r}{r_{i,\alpha}^\CR-\ee s}
    \leq \frac{\CR(q_{i,\alpha}^\CR)+\ee r_{i,\alpha}^\CR}{r_{i,\alpha}^\CR}\cdot\Big(1-\frac{\ee}{1-\alpha}\Big)^{-1},
\]
since $(1-\alpha)s \leq r_{i,\alpha}^\CR$. By taking the infimum over $i$, we obtain
\[
    \kappa_\beta(\CP) \leq \big(\kappa_\alpha(\CR)+\ee\big)\cdot\Big(1-\frac{\ee}{1-\alpha}\Big)^{-1},
\]
and by letting $\ee$ tend to $0$ we prove \eqref{Eq kappa essentiellement décroissante}. By symmetry, we also have $\kap_\beta(\CR) \leq \kap_\alpha(\CP)$, which yields $\kap_\alpha(\CP) = \kap_\alpha(\CR)$ at each point where $\alpha\mapsto \kappa_\alpha(\CR)$ is continuous. By letting successively $\beta$ and $\alpha$ tend to $\po$ in \eqref{Eq kappa essentiellement décroissante}, it follows that $\kap(\CP)\leq \kap(\CR)$. By symmetry we have $\kap(\CR)\leq \kap(\CP)$ and therefore $\kap(\CP) = \kap(\CR)$.
\end{Dem}

\begin{Prop}
\label{CEx Prop kap = beta}
Let $\bu=(1,\xi,\eta)$ where $\xi,\eta$ are non-zero real numbers and define the functions $\CL_i$, $\CL_i^*$ ($i=1,2,3$) as in Section~\ref{subsection setting} (with respect to $\bu$). Then
\begin{align}
\label{CEx Eq kap_0 = beta_0}
    \kap(\CL_3)=-\kap^*(\CL_1^*) = \frac{\lambdaL(\xi,\eta)}{1+\lambdaL(\xi,\eta)}.
\end{align}
\end{Prop}

\begin{Dem} Mahler's duality implies that $|\CL_3+\CL_1^*|$ is bounded. By Lemma~\ref{Lemma kappa(P) = kappa(R)} we conclude that $ \kap(\CL_3) = -\kap^*(\CL_1^*)$. Now, let us define $\NNN$ and $\Delta$ by $\NNN(\bx) = \norm{\bx}$ and $\Delta = \norm{\bu\wedge\bx}$ ($\bx\in\RR^3$). For each $\eee\geq 0$ with $\eee<\lambda(\xi,\eta)$, we denote by $\hnu_\eee(\xi,\eta)$ the exponent associated to $\NNN$ and $\Delta$ as in Section~\ref{section exposants lambda epsilon}. Let $(\bx_i)_{i\geq 0}$ be a sequence of minimal points associated to $\NNN$ and $\Delta$ and let us write $X_i:=\norm{\bx_i}$, $\Delta_i:=\norm{\bx_i\wedge\bu}$ ($i\geq 0$). A wellknown result in parametric geometry of numbers (see \cite[§4]{Schmidt2009}) states that:
\begin{align}
\label{Eq inter 1 dem relation kappa nu de espilon}
    \CL_1^*(q) = \min_{i\in\NN} \CL_{\bx_i}^*(q) \quad (q>0),
\end{align}
where $\CL_{\bx_i}^*(q) = \max\big(\log\Delta_i,\log X_i -q \big)$ (see \eqref{Eq formule pour L_bx^*}). Let us fix $0\leq \alpha < \po(-\CL_1^*) = -\pu_1^*$ and set $\eee:=\alpha/(1-\alpha)\geq 0$. Let us first prove that
\begin{align}
\label{Eq relation kappa alpha et nu epsilon}
    \kappa_{\alpha}(-\CL_1^*) = \frac{\hnu_\eee(\xi,\eta)}{1+\hnu_\eee(\xi,\eta)}.
\end{align}
We denote by $(q_i)_{i\geq 0}$, $(q_{k,\alpha})_k$ and $(r_{k,\alpha})_k$ the sequences associated to $\kappa_{\alpha}(-\CL_1^*)$ by Definition~\ref{Def kappa(P)}. Eq.~\eqref{Eq inter 1 dem relation kappa nu de espilon} implies that $q_i$ is the point at which $\CL_{\bx_i}^*$ changes slope (from $-1$ to $0$), which is precisely $\log(X_i)-\log(\Delta_i)$.
Let $i_1<i_2<\dots$ denote the sequence of indices $i$ such that $\Delta_i \leq \CX_i^{-\eee}$.
We claim that the sequence $(q_{k,\alpha})_k$ is the sequence $(q_{i_k})_k$. Indeed, the condition $-\CL_1^*(q_k)/q_k \geq \alpha$ is equivalent to the condition $\Delta_k \leq X_k^{-\eee}$, by using $-\CL_1^*(q_k) = -\log(\Delta_k)$. This implies that $r_{k,\alpha} = \log(X_{i_{k+1}})-\log(\Delta_{i_k})$, and we thus have
\begin{align}
\label{Eq inter 3 dem relation kappa nu de espilon}
    \kappa_\alpha(-\CL_1^*) = \liminf_{k\rightarrow\infty} \frac{-\CL_1^*(q_{i_k})}{r_{k,\alpha}} = \liminf_{k\rightarrow\infty} \frac{-\log(\Delta_{i_k})}{\log(X_{i_{k+1}})-\log(\Delta_{i_k})}.
\end{align}
Eqns.~\eqref{Eq formule nu epsilon comme lim inf avec log} and \eqref{Eq inter 3 dem relation kappa nu de espilon} together give \eqref{Eq relation kappa alpha et nu epsilon}. We conclude by noticing that when $\alpha$ tends to $-\pu_1^* = \po_3$, then $\eee$ tends to $-\pu_1^*/(1+\pu_1^*) = \lambda(\xi,\eta)$ by \eqref{eq dico} and \eqref{Eq dualité exposants}.
\end{Dem}

\section{Proof of Theorem \ref{Thm general spectre (lambda_min,lambda)}}
\label{section proof thm principal}

Recall that the first part of Theorem \ref{Thm general spectre (lambda_min,lambda)} follows from Laurent's inequalities \eqref{Eq relation Laurent} and from the fact that if $\lambda(\xi,\eta) = 1/2$, then $\lambdaL(\xi,\eta) = 1/2$. Let us prove the second part of Theorem \ref{Thm general spectre (lambda_min,lambda)}. Theorem~\ref{Thm valeur générique lambda_min} implies that there exist real numbers $\xi,\eta$ with $1,\xi,\eta$ linearly independent over $\QQ$, such that $\lambda(\xi,\eta) = \lambdaL(\xi,\eta) = 1/2$. Now, let $\hlambda\in\RR$ and $\lambda\in\RR\cup\{+\infty\}$ satisfying \eqref{Eq Thm ppal}. The strategy of the proof is to construct a $3$-system $\bP = (\CP_1,\CP_2,\CP_3)$ such that
\begin{align}
\label{Eq inter -1 Thm ppal}
    \lim_{q\rightarrow\infty} \CP_1(q) = +\infty,\quad \po(\CP_3) = \frac{\lambda}{1+\lambda}\quad \textrm{and} \quad \kappa(\CP_3) = \frac{\hlambda}{1+\hlambda},
\end{align}
with the convention that $\lambda/(1+\lambda) = 1$ if $\lambda=+\infty$. If $\bP$ is as above, Theorem~\ref{Thm Roy conjecture S&S} gives a non-zero vector $\bu\in\RR^3$ such that $\norm{\bL_\bu - \bP}$ is bounded. Moreover, we may suppose $\bu = (1,\xi,\eta)$ with $1,\xi,\eta$ linearly independent over $\QQ$, since $\CP_1$ is not bounded. 
Then, Lemma~\ref{Lemma kappa(P) = kappa(R)}, Proposition~\ref{CEx Eq kap_0 = beta_0} and relation \eqref{eq dico} imply that $\lambda(\xi,\eta) = \lambda$ and $\lambdaL(\xi,\eta) = \hlambda$.\\
In order to cover the full joint spectrum of $(\lambda,\lambdaL)$ we distinguish between two cases. Our first construction deals with the case $\max(1-\lambda,0) < \hlambda$ (note that this inequality is fulfilled if $\hlambda\geq 1/2$, since $\lambda>1/2$) and the second one deals with the case $\hlambda \leq 1/2$.\\

\noindent\textbf{First case.} Suppose that $\lambda,\hlambda$ satisfy $1<\lambda+\hlambda$ and $0<\hlambda$. For convenience, let us define $\nu\in(0,1/\hlambda]$ by $1/\nu = \hlambda\big(1+1/\lambda\big)\big(1+\hlambda/\lambda\big)$. This number satisfies the relations
\begin{align}
\label{Eq inter thm ppal calcul psi via nu}
1 - \frac{\hlambda}{\lambda}\nu - \Big(\frac{\hlambda}{\lambda}\Big)^2\nu = \frac{\lambda}{1+\lambda} \quad \textrm{and} \quad 1 + \nu - \Big(\frac{\hlambda}{\lambda}\Big)^2\nu = \frac{\lambda}{1+\lambda}\cdot \frac{1+\hlambda}{\hlambda}.
\end{align}
Under our hypotheses on $\lambda$ and $\hlambda$ we have
\begin{align}
\label{Eq inter 0 Thm ppal}
\frac{\hlambda}{\lambda} < \frac{1}{\nu} - \frac{\hlambda}{\lambda} - \Big(\frac{\hlambda}{\lambda}\Big)^2 \leq 1.
\end{align}
Indeed, the first inequality of \eqref{Eq inter 0 Thm ppal} is equivalent to $1 < \lambda + \hlambda$ and the second one is equivalent to the third inequality of \eqref{Eq Thm ppal}. Let $(\beta_k)_{k\geq 0}$ be a non-decreasing sequence of real numbers $>1$ such that $\beta_k$ tends to $\lambda/\hlambda\in(1,+\infty]$ as $k$ tends to infinity.
If $\lambda = +\infty$, we may take $\beta_k = k+1$ for each $k\geq 0$. If $\lambda<+\infty$, then we may simply take $\beta_k = \lambda/\hlambda$ for each $k\geq 0$. Since $\lambda/\hlambda>1$, the sequence $(q_k)_{k\geq 0}$ defined by $q_k := \prod_{i= 0}^{k}\beta_i$ tends to infinity. By \eqref{Eq inter 0 Thm ppal} and by the choice of $(\beta_k)_k$, there is an index $N\geq 1$ such that for each $k\geq N$ we have:
\begin{align}
\label{Eq inter 1 proof thm ppal}
    \frac{1}{\beta_k\beta_{k-1}} \leq \frac{1}{\beta_{k}} < \frac{1}{\nu} - \frac{1}{\beta_{k}} - \frac{1}{\beta_{k}\beta_{k-1}} \leq 1.
\end{align}

\noindent For each $k\geq 1$, we define a point $\ba^{(k)}=(a_1^{(k)},a_2^{(k)},a_3^{(k)})\in\RR^3$  by
\[
    \ba^{(k)} = q_k\times\Big(\frac{\nu}{\beta_k\beta_{k-1}},\frac{\nu}{\beta_{k}}, 1 - \frac{\nu}{\beta_{k}} - \frac{\nu}{\beta_k\beta_{k-1}}\Big).
\]
Note that $a_1^{(k+1)} = a_2^{(k)}$ since $q_{k+1} = \beta_{k+1} q_k$, and that $a_1^{(k)}+a_2^{(k)}+a_3^{(k)} = q_k$. Inequalities \eqref{Eq inter 1 proof thm ppal} may be rewritten as $a_1^{(k)}\leq a_2^{(k)} < a_3^{(k)} \leq a_2^{(k+1)}$ for each $k\geq N$. We now construct the $3$-system $\bP$ on $[q_N,+\infty)$ whose combined graph is shown on figure~\ref{figure 3-system P_alpha}.

\tikzstyle{alphaDte}=[dashed, thick,black] 
\tikzstyle{AutreDte}=[dashed, very thin, black!50] 
\tikzstyle{DteBas}=[dashed, very thin, black] 
\tikzstyle{Abscisses}=[dotted] 
\tikzstyle{ProlongationKappa}=[densely dotted] 
\tikzstyle{ProlongationKappaGros}=[dotted] 

\newcommand\PasVague{0.1}
\newcommand\VagueHautONE{-- ++ (\PasVague,\PasVague) -- ++(\PasVague,0)}
\newcommand\VagueHautFIVE{\VagueHautONE \VagueHautONE \VagueHautONE \VagueHautONE \VagueHautONE}
\newcommand\VagueHautTEN{\VagueHautFIVE \VagueHautFIVE}

\newcommand\VagueBasONE{-- ++ (\PasVague,0) -- ++(\PasVague,\PasVague)}
\newcommand\VagueBasFIVE{\VagueBasONE \VagueBasONE \VagueBasONE \VagueBasONE \VagueBasONE}
\newcommand\VagueBasTEN{\VagueBasFIVE \VagueBasFIVE}

\begin{figure}[H]
\begin{tabular}{lr}
    \begin{tikzpicture}[scale=0.6]
         \clip (7,-1) rectangle (17.8,9);


        \draw (7.8,1+1/3) -- (8,1+1/3) node [gray] {$\bullet$} -- (9+1/3,2+2/3)-- (16,2+2/3) node [gray] {$\bullet$} -- ++(0.3,0.3);
        \draw (7.5,1+1/3+0.3) node [] {$a_1^{(k)}$}; 
        \draw (17,2+2/3+0.3) node [] {$a_1^{(k+1)}$}; 


        \draw (7.8,2+2/3) -- (8,2+2/3) node [gray] {$\bullet$} -- (9+1/3,2+2/3) -- (10+2/3,4);
        \draw (7.5,2+2/3+0.3) node [] {$a_2^{(k)}$}; 
        \draw (17,5+1/3+0.3) node [] {$a_2^{(k+1)}$}; 

        \draw [thin] (10+2/3,4) \VagueBasTEN \VagueBasONE \VagueBasONE \VagueBasONE -- ++ (0.033057,0.033057) -- ++(0.033057,0);
        \draw [thin] (10+2/3,4) \VagueHautTEN \VagueHautONE \VagueHautONE \VagueHautONE -- ++ (0.033057,0.033057) -- ++(0.033057,0);

        \draw (13+1/3,5+1/3) -- (16,5+1/3) node [gray] {$\bullet$} -- ++(0.3,0);


        \draw (8,4) -- ++ (-0.3,-0.3);
        \draw (7.5,4+0.3) node [] {$a_3^{(k)}$}; 
        \draw (17,8+0.3) node [] {$a_3^{(k+1)}$}; 

        \draw (8,4) node [gray] {$\bullet$} -- (10+2/3,4);
        \draw (13+1/3,5+1/3)-- (16,8) node [gray] {$\bullet$} -- ++(0.3,0); 


        \draw [ProlongationKappaGros] (10+2/3,4)-- (12,4) -- (13+1/3,5+1/3); 
        \draw [ProlongationKappaGros] (10+2/3,4)-- (12,5+1/3) -- (13+1/3,5+1/3);

        \draw[Abscisses]  (8,4)--(8,0)node[below]{$q_{k}$};
        \draw[Abscisses]  (16,8)--(16,0)node[below]{$q_{k+1}=\beta_{k+1}q_k$};
        \draw[Abscisses] (6,2) --(6,0) node[below] {$r_{k-1}$};
        \draw[Abscisses] (12,4) node [] {$\bullet$} --(12,0) node[below] {$r_{k}$};

        \draw[Abscisses] (10+2/3,4)--(10+2/3,0) node[below] {$s_k$};
        \draw[Abscisses] (13+1/3,5+1/3)--(13+1/3,0) node[below] {$t_k$};


    \end{tikzpicture}

    &


    \begin{tikzpicture}[scale=0.35]
         \clip (2,-2) rectangle (22,14);

        \draw (1,1/6)-- (1+1/6,1/3)-- (2,1/3)-- (2+1/3,2/3)-- (4,2/3)-- (4+2/3,1+1/3)-- (8,1+1/3)-- (9+1/3,2+2/3)-- (16,2+2/3)-- (18+2/3,5+1/3)-- (24,5+1/3) -- ++ (0.5,0);

        \draw (1,1/3)-- (1+1/6,1/3)-- (1+1/3,0.5)-- (1.5,0.5)  -- (1+2/3,2/3)-- (2+1/3,2/3)-- (2+2/3,1)-- (2.779770672705434,1)-- (2.8961792298261355,1.115955888587084)-- (3.005768087333419,1.115955888587084)-- (3.1123903876198655,1.2221635753647355)-- (3.2221635753647355,1.2221635753647355)-- (3+1/3,1+1/3)-- (4+2/3,1+1/3) -- (5+1/3,2);

        \draw [thin] (5+1/3,2) \VagueBasFIVE \VagueBasONE -- ++ (0.06717,0.06717) -- ++(0.06717,0);
        \draw [thin] (5+1/3,2) \VagueHautFIVE \VagueHautONE -- ++ (0.06717,0.06717) -- ++(0.06717,0);

        \draw (6+2/3,2+2/3) -- (9+1/3,2+2/3) -- (10+2/3,4);

        \draw [thin] (10+2/3,4) \VagueBasTEN \VagueBasONE \VagueBasONE \VagueBasONE -- ++ (0.033057,0.033057) -- ++(0.033057,0);
        \draw [thin] (10+2/3,4) \VagueHautTEN \VagueHautONE \VagueHautONE \VagueHautONE -- ++ (0.033057,0.033057) -- ++(0.033057,0);

        \draw (13+1/3,5+1/3) -- (18+2/3,5+1/3) -- (21+1/3,8);

        \draw [thin] (21.321490687630142,8) \VagueHautTEN \VagueHautFIVE \VagueHautONE ;
        \draw [thin] (21.321490687630142,8) \VagueBasTEN \VagueBasFIVE \VagueBasONE ;

        \draw (1,0.5)-- (1+1/3,0.5)-- (1.5,2/3)-- (1+2/3,2/3)-- (2,1)-- (2+2/3,1)-- (2.7749928750570105,1.115955888587084)-- (2.8961792298261355,1.115955888587084)-- (3.002801530112582,1.2221635753647355)-- (3.1123903876198655,1.2221635753647355)-- (3.217991997277099,1.327354540590206) -- (3.3281489250350873,1.3281489250350875);

        \draw (3.3281489250350873,1.3281489250350875)-- (4,2) node [gray] {$$} -- (5.334905988262791,2);
        \draw (6.669244470632253,2.6692444706322536)-- (8,4) node [gray] {$$} -- (10.667056956092855,4);
        \draw (13.333170823451885,5.333170823451886)-- (16,8) node [gray] {$$} -- (21.321490687630142,8);

        \draw [ProlongationKappa] (2+2/3,1)-- (3,1) -- (3+1/3,1+1/3);
        \draw [ProlongationKappa] (5+1/3,2)-- (6,2) node {$\bullet$} -- (6+2/3,2+2/3);
        \draw [ProlongationKappa] (10+2/3,4)-- (12,4) node {$\bullet$} -- (13+1/3,5+1/3);
        \draw [ProlongationKappa] (21.321490687630142,8)-- (24,8)  node {$\bullet$} -- ++ (0.5,0.5);

        \draw[Abscisses]  (8.00015103975756,4)--(8.00015103975756,0)node[below]{$q_{k}$};
        \draw[Abscisses]  (15.983524939035563,8)--(15.983524939035563,0)node[below]{$q_{k+1}$};
        \draw[Abscisses] (6,2)--(6,0)node[below]{$r_{k-1}$};
        \draw[Abscisses] (12,4)--(12,0)node[below]{$r_{k+1}$};


    \end{tikzpicture}
    \end{tabular}

    \caption{\label{figure 3-system P_alpha}
    combined graph of a $3$-system $\bP$}
\end{figure}
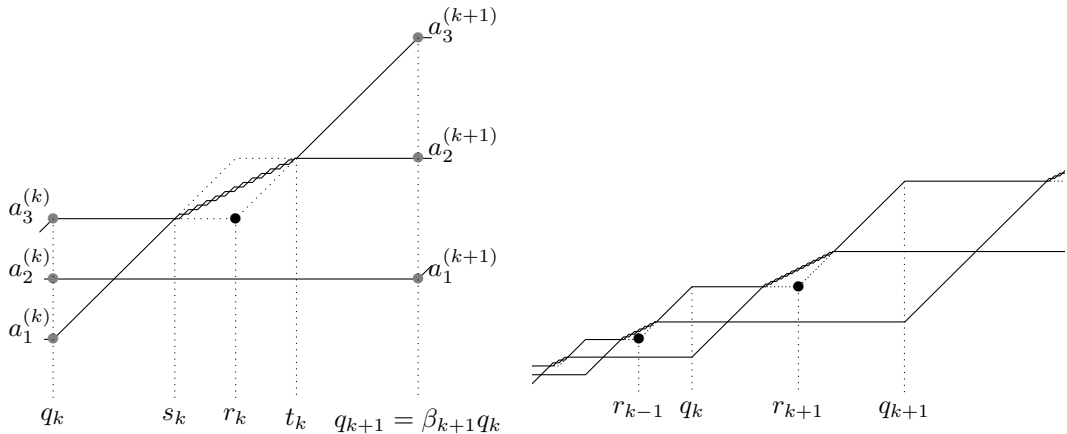

\noindent Set $\Delta=\{(x_1,x_2,x_3)\in\RR^3\;;\; x_1\leq x_2\leq x_3 \}$ and denote by $\Phi:\RR^3\rightarrow\Delta$ the continuous map which lists the coordinates of a point in monotone non-decreasing order. Let $s_k$ and $t_k$ be such that $a_1^{(k)} + s_k-q_{k} = a_3^{(k)}$ and $a_2^{(k+1)} = a_3^{(k+1)} - (q_{k+1} - t_k)$. We have
\[
    s_k = \Big(2-\frac{\nu}{\beta_k} - \frac{2\nu}{\beta_k\beta_{k-1}}\Big)q_k\quad \textrm{and}\quad t_k = \Big(2\nu+\frac{\nu}{\beta_k}\Big)q_k,
\]
and thus $s_k\leq t_k$ thanks to the last inequality of \eqref{Eq inter 1 proof thm ppal}. We define
\[
    \bP(q)
    =\left\{
    \begin{array}{ll}
    \Phi\big(a_1^{(k)}+q-q_{k},a_2^{(k)},a_3^{(k)}\big) & \textrm{if $q_{k}\leq q \leq s_k$}\\
    \Phi\big(a_1^{(k+1)},a_2^{(k+1)},a_3^{(k+1)}+q-q_{k+1}\big)& \textrm{if $t_k\leq q < q_{k+1}$}.
    \end{array}
    \right.
\]
In order to define $\bP$ on $[s_k,t_k]$, note that the ratio $a_2^{(k+1)}/t_k$ is smaller than $1/2$ and tends to $1/(2+\hlambda/\lambda)$ as $k$ tends to infinity, whereas the ratio $a_3^{(k+1)}/q_{k+1}$ tends to $\lambda/(1+\lambda)$ by using \eqref{Eq inter thm ppal calcul psi via nu}. Yet, the inequality $1<\lambda + \hlambda$ implies that the first limit is less than the second one. There exists therefore a real number $\theta$ such that
\[
    \frac{a_2^{(k+1)}}{t_k} < \theta < \frac{a_3^{(k+1)}}{q_{k+1}}
\]
for $k$ large enough. For each $q\in[s_k,t_k]$, we set $\CP_1(q) = a_2^{(k)}$ and we define $\CP_2(q)$ and $\CP_3(q)$ such that $\bP = (\CP_1,\CP_2,\CP_3)$ satisfies the hypotheses of a $3$-system (which is possible since the line passing through the points $(s_k,a_3^{(k)})$ and $(t_k,a_2^{(k+1)})$ has slope $1/2$) and such that, when $k$ is large enough, we have
\begin{align}
\label{Eq inter 3 proof thm ppal}
    \frac{\CP_3(q)}{q} < \theta \quad \textrm{for each $q\in [s_k,t_k]$}
\end{align}
(see figure~\ref{figure 3-system P_alpha}). Let $r_k$ be the abscissa of the intersection of the horizontal line passing through $(q_k,\CP_3(q_k))$ and of the line with slope $1$ passing through $(q_{k+1},\CP_3(q_{k+1}))$. We have
\[
    r_k = a_3^{(k)}-a_3^{(k+1)}+q_{k+1} = \Big(1+\nu-\frac{\nu}{\beta_k\beta_{k-1}}\Big)q_k.
\]
By \eqref{Eq inter 3 proof thm ppal} and \eqref{Eq inter thm ppal calcul psi via nu}, for each $\theta\leq \alpha < \lambda/(1+\lambda)$ we have
\begin{align*}
    \po(\CP_3) = \limsup_{k\rightarrow\infty} \frac{\CP_3(q_k)}{q_k} = \limsup_{k\rightarrow\infty} \frac{a_3^{(k)}}{q_k} = \frac{\lambda}{1+\lambda} \quad \textrm{and}\quad \kappa_\alpha(\CP_3) = \kappa(\CP_3) = \liminf_{k\rightarrow\infty} \frac{\CP_3(q_k)}{r_k} = \frac{\hlambda}{1+\hlambda}.
\end{align*}
Thus, $\bP$ satisfies \eqref{Eq inter -1 Thm ppal}, which concludes the first case. 

\bigskip

\noindent\textbf{Second case.} Suppose that $\hlambda\leq 1/2$. Under this additional condition, \eqref{Eq Thm ppal} may simply be rewritten as $0\leq \hlambda \leq \frac{1}{2} < \lambda \leq +\infty$, which is equivalent to
\[
    0\leq \frac{\hlambda}{1+\hlambda} \leq \frac{1}{3} < \frac{\lambda}{1+\lambda} \leq 1.
\]
Fix $\theta\in\RR$ such that $1/3<\theta<\lambda/(1+\lambda)$. Let $(\alpha_k)_{k\geq 1}$, $(\psi_k)_{k\geq 1}$ be two sequences of real numbers which tend to $\hlambda/(1+\hlambda)$ and $\lambda/(1+\lambda)$ respectively, and such that  for each $k\geq 1$, we have
\[
    0 < \alpha_k \leq \frac{1}{3} < \theta < \psi_k < 1.
\]
Let $(q_k)_{k\geq 0}$ be the sequence defined by $q_0 = 1$ and
\begin{align}
\label{Eq inter 10 Thm ppal}
q_{k+1} = \frac{\psi_k}{1-\psi_{k+1}}\Big(\frac{1}{\alpha_k} -1\Big)q_k\quad (k\geq 0).
\end{align}
Note that $q_{k+1}/q_k > 2\theta/(1-\theta) > 1$ for each $k\geq 0$, which implies that the sequence $(q_k)_k$ tends to infinity. For each $k\geq 0$, let us define the abscissas $s_k$ and $t_k$ by $s_k/3 = \psi_kq_k$ and $t_k/3 = (1-\psi_{k+1})q_{k+1}/2$. Let $r_k$ be the abscissa of the intersection point of the horizontal line passing through $(q_k,\psi_kq_k)$ and of the line with slope $1$ passing through $(q_{k+1},\psi_{k+1}q_{k+1})$ (see figure~\ref{figure 3-system P_alpha construction facile}). We have $r_k = \psi_kq_k + (1-\psi_{k+1})q_{k+1}$, which may be rewritten as $\alpha_kr_k = \psi_kq_k$ by \eqref{Eq inter 10 Thm ppal}. Since $0<\alpha_k \leq 1/3$, we have $s_k\leq r_k\leq t_k$. Now, let $\bP=(\CP_1,\CP_2,\CP_3)$ be a $3$-system on $[q_0,+\infty[$ such that for each $k\geq 0$, we have


\begin{align*}
    \frac{\bP(q_k)}{q_k} = \Big(\frac{1-\psi_k}{2}, \frac{1-\psi_k}{2}, \psi_k \Big), \quad
    \frac{\bP(s_k)}{s_k} = \frac{\bP(t_k)}{t_k} = \Big(\frac{1}{3}, \frac{1}{3}, \frac{1}{3} \Big),
\end{align*}
and such that
\begin{align}
\label{Eq inter 11 Thm ppal}
\frac{\bP(q)}{q}\leq \theta \quad \textrm{for each $q\in[s_k,t_k]$}.
\end{align}
An example of such $3$-system is represented on figure~\ref{figure 3-system P_alpha construction facile}.

\tikzstyle{PremierMin}=[]
\tikzstyle{DeuxiemeMin}=[]
\tikzstyle{TroisiemeMin}=[]

\tikzstyle{Double}=[thick] 
\tikzstyle{Droite}=[dotted] 
\tikzstyle{Bullet}=[gray] 
\tikzstyle{Abscisses}=[dotted] 

\begin{figure}[H]
    \begin{tabular}{l|r}
        \begin{tikzpicture}[scale=0.50]

            \clip(0.6,-1) rectangle (14.3,8);

            \draw (11.5,6.5) node {slope $\frac{\lambda}{1+\lambda}$};

            \draw (12.5,4.5) node {slope $\frac{1}{3}$};
            \draw  (11.5,9.5/6) node {slope $\frac{\hlambda}{1+\hlambda}$};

            \draw [Droite, domain=2.3:13.9] plot(\x,{(1*\x)/2});
            \draw [Droite, domain=2.3:13.9] plot(\x,{(1*\x)/3});
            \draw [Droite, domain=2.3:13.9] plot(\x,{(1*\x)/6}); 


            \draw [PremierMin] (2+2/3,2/3) -- ++(-0.2,0);
            \draw [PremierMin] (2+2/3,2/3) -- (3+1/3,2/3)-- (4,1+1/3);
            \draw [PremierMin] (10,3+1/3)-- (13+1/3,3+1/3) -- ++(0.4,0);

            \draw [PremierMin] (4,1+1/3)-- (4.522225509595202,1.3333333333333333)-- (4.783338264392803,1.5944460881309344)-- (5.228759511017056,1.5944460881309344)-- (5.451470134329184,1.8171567114430613)-- (5.925628235574356,1.8171567114430613)-- (6.162707286196943,2.0542357620656477)-- (6.657391712171345,2.0542357620656477)-- (6.904733925158546,2.3015779750528487)-- (7.360418334147413,2.3015779750528487)-- (7.588260538641846,2.529420179547282)-- (8.048050212576559,2.529420179547282)-- (8.277945049543916,2.759315016514639)-- (8.75931501651464,2.759315016514639)-- (9,3)-- (9+2/3,3)-- (10,3+1/3);


            \draw [TroisiemeMin] (2+2/3,1+1/3) -- ++(-0.2,-0.2);
            \draw [TroisiemeMin] (2+2/3,1+1/3) node[Bullet] {$\bullet$} node [above left] {$\psi_kq_k$} -- (4,1+1/3);
            \draw [TroisiemeMin] (10,3+1/3)-- (13+1/3,6+2/3) node[Bullet] {$\bullet$} -- ++(0.4,0);

            \draw [TroisiemeMin] (4,1+1/3)-- (4.261112754797601,1.5944460881309344)-- (4.783338264392803,1.5944460881309344)-- (5.006048887704931,1.8171567114430613)-- (5.451470134329184,1.8171567114430613)-- (5.6885491849517695,2.0542357620656477)-- (6.162707286196943,2.0542357620656477)-- (6.410049499184144,2.3015779750528487)-- (6.904733925158546,2.3015779750528487)-- (7.132576129652978,2.529420179547282)-- (7.588260538641846,2.529420179547282)-- (7.818155375609203,2.759315016514639)-- (8.277945049543916,2.759315016514639)-- (8.518630033029277,3)-- (9,3)-- (9+1/3,3+1/3)-- (10,3+1/3);


            \draw [DeuxiemeMin] (2+2/3,2/3) -- ++(-0.2,0);
            \draw [DeuxiemeMin] (2+2/3,2/3) node[Bullet] {$\bullet$} -- (3+1/3,1+1/3)-- (4,1+1/3);
            \draw (2.9,0.5) node [left] {$\frac{1-\psi_k}{2}q_k$};
            \draw [DeuxiemeMin] (10,3+1/3)-- (13+1/3,3+1/3) node[Bullet] {$\bullet$} -- ++(0.4,0.4);

            \draw [DeuxiemeMin] (4,1.3333333333333333)-- (4.261112754797601,1.3333333333333333)-- (4.522225509595201,1.5944460881309346)-- (5.006048887704931,1.5944460881309344)-- (5.228759511017056,1.8171567114430611)-- (5.68854918495177,1.8171567114430613)-- (5.925628235574357,2.0542357620656477)-- (6.410049499184144,2.0542357620656477)-- (6.657391712171344,2.3015779750528487)-- (7.132576129652978,2.3015779750528487)-- (7.360418334147411,2.529420179547282)-- (7.818155375609203,2.529420179547282)-- (8.04805021257656,2.759315016514639)-- (8.518630033029277,2.759315016514639)-- (8.75931501651464,3)-- (9+1/3,3)-- (9+2/3,3+1/3)-- (10,3+1/3);


            \draw [dashed] (4.522225509595202,1+1/3)-- (8,1+1/3)-- (9+2/3,3);


            \draw [Double] (2+2/3,2/3) -- ++(-0.2,0);
            \draw [Double] (2+2/3,2/3) -- (2+2/3,2/3);
            \draw [Double] (3+1/3,1+1/3)-- (4,1+1/3);
            \draw [Double] (10,3+1/3)-- (13+1/3,3.+1/3);
            \draw [Double] (16+2/3,6+2/3)-- (20,6+2/3);


            \draw [Abscisses] (2+2/3,2/3) node [Bullet] {$\bullet$} -- (2+2/3,0) node [below] {$q_k$};
            \draw [Abscisses] (4,1+1/3) node [Bullet] {$\bullet$} -- (4,0) node [below] {$s_k$};
            \draw [Abscisses] (8,1+1/3) node [Bullet] {$\bullet$} -- (8,0) node [below] {$r_{k}$};
            \draw [Abscisses] (10,3+1/3) node [Bullet] {$\bullet$} -- (10,0) node [below] {$t_{k}$};
            \draw [Abscisses] (13+1/3,3+1/3) node [Bullet] {$\bullet$} -- (13+1/3,0) node [below] {$q_{k+1}$};
        \end{tikzpicture}
        &
        \begin{tikzpicture}[scale=0.27]

            \clip(1.6,-1.5) rectangle (26,13);


            \draw (21,11)  node {slope $\frac{\lambda}{1+\lambda}$};
            \draw (22,19/3-0.5)  node {slope $\frac{1}{3}$};
            \draw (21,19/6-0.5) node {slope $\frac{\hlambda}{1+\hlambda}$};

            \draw [Droite, domain=2.3:13.9] plot(\x,{(1*\x)/2});
            \draw [Droite, domain=2.3:13.9] plot(\x,{(1*\x)/3});
            \draw [Droite, domain=2.3:13.9] plot(\x,{(1*\x)/6}); 

            \draw [Droite, domain=-1.5:25] plot(\x,{(1*\x)/2});
            \draw [Droite, domain=-1.5:25] plot(\x,{(1*\x)/3});
            \draw [Droite, domain=-1.5:25] plot(\x,{(1*\x)/6}); 


            \draw [PremierMin] (2,2/3)-- (3+1/3,2/3)-- (4,1+1/3);
            \draw [PremierMin] (10,3+1/3)-- (16+2/3,3+1/3)-- (20,6+2/3);

            \draw [PremierMin] (4,1+1/3)-- (4.522225509595202,1.3333333333333333)-- (4.783338264392803,1.5944460881309344)-- (5.228759511017056,1.5944460881309344)-- (5.451470134329184,1.8171567114430613)-- (5.925628235574356,1.8171567114430613)-- (6.162707286196943,2.0542357620656477)-- (6.657391712171345,2.0542357620656477)-- (6.904733925158546,2.3015779750528487)-- (7.360418334147413,2.3015779750528487)-- (7.588260538641846,2.529420179547282)-- (8.048050212576559,2.529420179547282)-- (8.277945049543916,2.759315016514639)-- (8.75931501651464,2.759315016514639)-- (9,3)-- (9+2/3,3)-- (10,3+1/3);

            \draw [PremierMin] (20,6+2/3)-- (20.8146781471295,6+2/3)-- (21.22201722069425,7.074005740231417)-- (22.12911141089019,7.074005740231417)-- (22.58265850598816,7.527552835329387)-- (23.543111177960334,7.527552835329386)-- (24.02333751394642,8.007779171315473);

            \draw [PremierMin] (0.9003862815856366,0.30012876052854554)-- (1.1691496546287252,0.30012876052854554)-- (1.3035313411502694,0.4345104470500898)-- (1.5449426895031322,0.4345104470500898)-- (1.6656483636795634,0.5552161212265211)-- (1.8885494545598545,0.5552161212265211)-- (2,2/3);


            \draw [TroisiemeMin] (2,2/3)-- (2+2/3,1+1/3) node [Bullet] {$\bullet$} -- (4,1+1/3);
            \draw [TroisiemeMin] (10,3+1/3)-- (13+1/3,6+2/3) node [Bullet] {$\bullet$} -- (20,6+2/3);

            \draw [TroisiemeMin] (4,1+1/3)-- (4.261112754797601,1.5944460881309344)-- (4.783338264392803,1.5944460881309344)-- (5.006048887704931,1.8171567114430613)-- (5.451470134329184,1.8171567114430613)-- (5.6885491849517695,2.0542357620656477)-- (6.162707286196943,2.0542357620656477)-- (6.410049499184144,2.3015779750528487)-- (6.904733925158546,2.3015779750528487)-- (7.132576129652978,2.529420179547282)-- (7.588260538641846,2.529420179547282)-- (7.818155375609203,2.759315016514639)-- (8.277945049543916,2.759315016514639)-- (8.518630033029277,3)-- (9,3)-- (9+1/3,3+1/3)-- (10,3+1/3);

            \draw [TroisiemeMin] (20,6+2/3)-- (20.40733907356475,7.074005740231416)-- (21.22201722069425,7.074005740231417)-- (21.67556431579222,7.527552835329386)-- (22.58265850598816,7.527552835329387)-- (23.06288484197425,8.007779171315473)-- (24.02333751394642,8.007779171315473);

            \draw [TroisiemeMin] (0.9003862815856366,0.30012876052854554)-- (1.034767968107181,0.4345104470500898)-- (1.3035313411502694,0.4345104470500898)-- (1.4242370153267006,0.5552161212265211)-- (1.6656483636795634,0.5552161212265211)-- (1.7770989091197091,2/3)-- (2,2/3);


            \draw [DeuxiemeMin] (2,2/3)-- (2+2/3,2/3)-- (3+1/3,1+1/3)-- (4,1+1/3);
            \draw [DeuxiemeMin] (10,3+1/3)-- (13+1/3,3.+1/3)-- (16+2/3,6+2/3)-- (20,6+2/3);

            \draw [DeuxiemeMin] (4,1.3333333333333333)-- (4.261112754797601,1.3333333333333333)-- (4.522225509595201,1.5944460881309346)-- (5.006048887704931,1.5944460881309344)-- (5.228759511017056,1.8171567114430611)-- (5.68854918495177,1.8171567114430613)-- (5.925628235574357,2.0542357620656477)-- (6.410049499184144,2.0542357620656477)-- (6.657391712171344,2.3015779750528487)-- (7.132576129652978,2.3015779750528487)-- (7.360418334147411,2.529420179547282)-- (7.818155375609203,2.529420179547282)-- (8.04805021257656,2.759315016514639)-- (8.518630033029277,2.759315016514639)-- (8.75931501651464,3)-- (9+1/3,3)-- (9+2/3,3+1/3)-- (10,3+1/3);

            \draw [DeuxiemeMin] (20,6+2/3)-- (20.407339073564753,6.666666666666667)-- (20.814678147129502,7.074005740231416)-- (21.67556431579222,7.074005740231417)-- (22.129111410890182,7.527552835329387)-- (23.06288484197425,7.527552835329386)-- (23.543111177960338,8.007779171315473)-- (24.02333751394642,8.007779171315473);

            \draw [DeuxiemeMin] (0.9003862815856366,0.30012876052854554)-- (1.034767968107181,0.30012876052854554)-- (1.1691496546287254,0.4345104470500898)-- (1.4242370153267008,0.43451044705008984)-- (1.544942689503132,0.5552161212265211)-- (1.7770989091197096,0.5552161212265211)-- (1.8885494545598553,2/3)-- (2,2/3);


            \draw [dashed] (4.522225509595202,1+1/3)-- (8,1+1/3) node [Bullet] {$\bullet$} -- (9+2/3,3);


            \draw [Double] (2,2/3)-- (2+2/3,2/3);
            \draw [Double] (3+1/3,1+1/3)-- (4,1+1/3);
            \draw [Double] (10,3+1/3)-- (13+1/3,3.+1/3);
            \draw [Double] (16+2/3,6+2/3)-- (20,6+2/3);

            \draw [Abscisses] (2+2/3,2/3) node [Bullet] {$\bullet$} -- (2+2/3,0) node [below] {$q_k$};
            \draw [Abscisses] (8,1+1/3) node [Bullet] {$\bullet$} -- (8,0) node [below] {$r_{k}$};
            \draw [Abscisses] (13+1/3,3+1/3) node [Bullet] {$\bullet$} -- (13+1/3,0) node [below] {$q_{k+1}$};

        \end{tikzpicture}
    \end{tabular}

    \caption{\label{figure 3-system P_alpha construction facile}
     combined graph of the $3$-system $\bP$}
\end{figure}
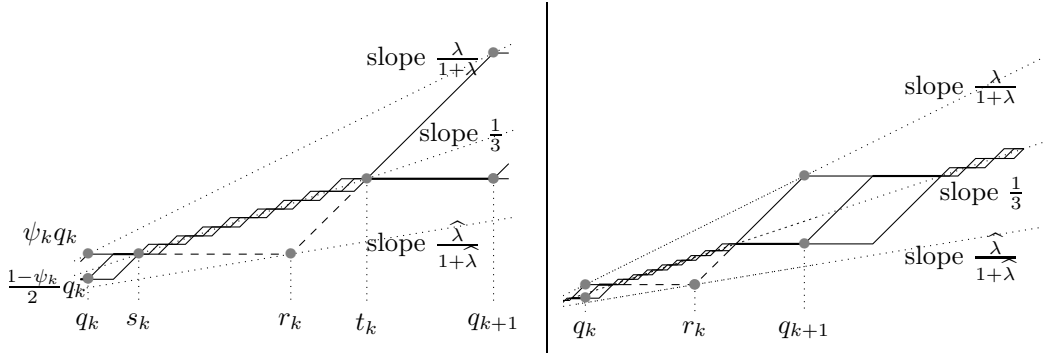

By \eqref{Eq inter 11 Thm ppal} and since $\theta \leq \psi_k$ for each $k\geq 0$, it is clear that such a $3$-system $\bP$ satisfies $\po(\CP_3) = \lambda/(1+\lambda)$. Moreover, \eqref{Eq inter 11 Thm ppal} also implies that $\kappa_{\alpha}(\CP_3) = \hlambda/(1+\hlambda)$ for each $\alpha$ such that $\theta<\alpha < \lambda/(1+\lambda)$. We thus have $\kappa(\CP_3) = \hlambda/(1+\hlambda)$ and $\bP$ satisfies \eqref{Eq inter -1 Thm ppal}. 
This ends the proof of Theorem~\ref{Thm general spectre (lambda_min,lambda)}.

\section{Final remarks and open questions}
\label{section: finals remarks and open questions}

We are grateful to the referee who suggested most of the following remarks and questions.\\

As mentioned in the introduction, it would be interesting to study the spectrum of $(\lambda,\hlambda,\lambdaL)$, or to begin with that of $(\hlambda,\lambdaL)$.

Note that the spectrum of $\hlambda-\lambdaL$ is simpler to study. Since $0\leq \lambdaL \leq \hlambda \leq 1$, it is included in $[0,1]$. By using the first construction from Section~\ref{section proof thm principal} with $\lambda = +\infty$ and $0<\hlambda \leq 1$, we obtain a $3$-system $(\CP_1,\CP_2,\CP_3)$ with $\pu(\CP_3) = 1/2$ and $\kap(\CP_3) = \hlambda/(1+\hlambda)$. It corresponds to a pair $(\xi,\eta)$ with $\xi,\eta,1$ linearly independent over $\QQ$, such that $\hlambda(\xi,\eta) = 1$ and $\lambdaL(\xi,\eta) = \hlambda$. Thus the spectrum of $\hlambda-\lambdaL$ contains $[0,1)$. It seems to be possible to modify the first construction to obtain $\hlambda(\xi,\eta) = 1$ and $\lambdaL(\xi,\eta) = 0$, which would show that the spectrum of $\hlambda-\lambdaL$ is equal to the full interval $[0,1]$.

The exponents $\hlambda_\mu$, $\lambdaL$ can be easily defined for $n$ numbers $\xi_1,\dots,\xi_n$ ($n\geq 2$). What is the spectrum of $(\lambda,\lambdaL)$ in this general setting?

What is the Hausdorff dimension of the sets $\{(\xi,\eta) \; |\; \lambdaL(\xi,\eta) = \hlambda\}$ (for $0\leq \hlambda \leq 1$)? We could maybe provide an answer to this question by using methods used in \cite{das2019variational}.

\bibliographystyle{abbrv}


\begin{thebibliography}{10}

\bibitem{bugeaud2015exponents}
Y.~Bugeaud.
\newblock Exponents of {D}iophantine approximation.
\newblock In D.~Badziahin, A.~Gorodnik, and N.~Peyerimhoff, editors, {\em
  Dynamics and Analytic Number Theory}, volume 437 of {\em London Math. Soc.
  Lecture Note Ser.}, chapter~2, pages 96--135. Cambridge University Press,
  2016.

\bibitem{bugeaud2005exponentsSturmian}
Y.~Bugeaud and M.~Laurent.
\newblock Exponents of {D}iophantine approximation and {S}turmian continued
  fractions.
\newblock {\em Ann. Inst. Fourier}, 55(3):773--804, 2005.

\bibitem{bugeaud2005exponents}
Y.~Bugeaud and M.~Laurent.
\newblock On exponents of homogeneous and inhomogeneous {D}iophantine
  approximation.
\newblock {\em Moscow Math. J}, 5(4):747--766, 2005.

\bibitem{cassaigne1999limit}
J.~Cassaigne.
\newblock Limit values of the recurrence quotient of {S}turmian sequences.
\newblock {\em Theoret. Comput. Sci.}, 218(1):3--12, 1999.

\bibitem{das2019variational}
T.~Das, L.~Fishman, D.~Simmons, and M.~Urba{\'n}ski.
\newblock A variational principle in the parametric geometry of numbers.
\newblock {\em arXiv preprint arXiv:1901.06602v3}, 2019.

\bibitem{davenport1967approximation}
H.~Davenport and W.~Schmidt.
\newblock Approximation to real numbers by quadratic irrationals.
\newblock {\em Acta Arith.}, 13(2):169--176, 1967.

\bibitem{davenport1969approximation}
H.~Davenport and W.~Schmidt.
\newblock Approximation to real numbers by algebraic integers.
\newblock {\em Acta Arith.}, 15(4):393--416, 1969.

\bibitem{fischler2007palindromic}
S.~Fischler.
\newblock Palindromic prefixes and diophantine approximation.
\newblock {\em Monatsh. Math.}, 151(1):11--37, 2007.

\bibitem{laurent2006exponents}
M.~Laurent.
\newblock Exponents of {D}iophantine approximation in dimension two.
\newblock {\em Canad. J. Math}, 61:165--189, 2009.

\bibitem{poels2017exponents}
A.~Po{\"e}ls.
\newblock Exponents of {D}iophantine approximation in dimension 2 for numbers
  of {S}turmian type.
\newblock {\em Math. Z.}, https://doi.org/10.1007/s00209-019-02280-2, 2019.

\bibitem{poelsPhD}
A.~Po{\"e}ls.
\newblock {\em Applications de la géométrie paramétrique des nombres à
  l'approximation diophantienne}.
\newblock PhD thesis, Université {P}aris-Sud,
  https://tel.archives-ouvertes.fr/tel-01827304/document, 2018.

\bibitem{Roy_juin}
D.~Roy.
\newblock {On Schmidt and Summerer parametric geometry of numbers}.
\newblock {\em Ann. of Math.}, 182:739--786, 2015.

\bibitem{Roy_octobre}
D.~Roy.
\newblock Spectrum of the exponents of best rational approximation.
\newblock {\em Math. Z.}, 283(1-2):143--155, 2016.

\bibitem{schleischitz2016spectrum}
J.~Schleischitz.
\newblock On the spectrum of {D}iophantine approximation constants.
\newblock {\em Mathematika}, 62(1):79--100, 2016.

\bibitem{Schmidt2009}
W.~M. Schmidt and L.~Summerer.
\newblock Parametric geometry of numbers and applications.
\newblock {\em Acta Arith.}, 140:67--91, 2009.

\bibitem{Schmidt2013}
W.~M. Schmidt and L.~Summerer.
\newblock Diophantine approximation and parametric geometry of numbers.
\newblock {\em Monatsh. Math.}, 169:51--104, 2013.

\end{thebibliography}
\footnotesize {

}

\Addresses

\end{document}